\documentclass[a4paper,leqno,final]{amsart}
\usepackage{amsfonts,epic,eepic}
\usepackage{palatino}
\usepackage{amsmath}
\usepackage{a4wide}
       
\theoremstyle{plain}
\begingroup

\newtheorem{teo}{Theorem}[section]
\newtheorem{lemma}[teo]{Lemma}
\newtheorem{prop}[teo]{Proposition}
\newtheorem{cor}[teo]{Corollary}
\endgroup

\theoremstyle{definition}

\newtheorem{dfnz}[teo]{Definition}
\newtheorem{conge}[teo]{Conjecture}

\theoremstyle{remark}

\newtheorem{rem}[teo]{Remark}

\numberwithin{equation}{section}

\def\R{{{\mathbb R}}}
\def\SS{{{\mathbb S}}}

\def\H{{{\mathcal H}}}
\def\NN{{{\mathbb N}}}
 
\def\cutt{{Cut(K)}}
\def\lra{\to}
\def\Derpar#1 { \frac{\partial~ } {\partial {#1} }} 
\def\derpar#1 #2 { \frac{\partial{#2}} {\partial {#1} }}

\def\real{\R}

\def\rl{\real}

\def\overbar{\overline}

\def\infinity{\infty}
\def\supdiff{{\partial^+}}
\def\subdiff{{\partial^-}}

\def\bordreg{K}
\def\uzero{{u_0}}
\def\composed{\circ}
\def\chart{\psi}
\def\l{\lambda}

\begin{document}

\title[Hamilton--Jacobi Equations and Distance
Functions]{Hamilton--Jacobi Equations and Distance Functions on
  Riemannian Manifolds}
\author[Carlo Mantegazza]{Carlo Mantegazza}
\address[Carlo Mantegazza]{Scuola Normale Superiore Pisa, Italy, 56126}
\email[C. Mantegazza]{mantegaz@sns.it}
\author[Andrea Carlo Mennucci]{Andrea Carlo Mennucci}
\address[Andrea Mennucci]{Scuola Normale Superiore Pisa, Italy, 56126}
\email[A. C. Mennucci]{mennucci@sns.it}
\date{\today}

\begin{abstract} The paper is concerned with the 
properties of the distance function from a closed subset of a
Riemannian manifold, with particular attention to the set of
singularities.
\end{abstract}

\maketitle

\tableofcontents

\section{Introduction}

We are concerned with the properties of the singular set of the
distance function from a closed subset of a $n$--dimensional, 
smooth and  connected Riemannian manifold $(M,g)$, with particular
attention to its rectifiability.

\begin{dfnz}\label{rectdef} We say that a subset $S$ of $(M,g)$ is 
$C^r$--{\em rectifiable}, with $r\geq1$, if it can be covered by a
countable family of embedded $C^r$ submanifolds of dimension $(n-1)$,
with the exception of a set of $\H^{n-1}$ zero measure, where
$\H^{n-1}$ is the $(n-1)$--dimensional Hausdorff measure on $M$.\\
We will simply say that a set is {\em rectifiable} when it is at least
$C^1$--{\em rectifiable}.
\end{dfnz}
See~\cite{fede,simon} for a complete discussion of the notion of
  rectifiability.

The distance function from a closed, not empty subset $K$ of $(M,g)$ 
is defined in the usual way,
\begin{equation*}
d_K(x)=\inf_{y\in K}d(x,y)
\end{equation*}
where $d$ is the distance on $M$ induced by the metric tensor $g$.\\
The  singular set $Sing$ of $d_K:M\to\R$ is the set where this
function fails to be differentiable.\\
Our study of the rectifiability of $Sing$ relies on the theory of {\em viscosity
solutions} of Hamilton--Jacobi equations. Indeed, we show that the distance
function $d_K$ is a viscosity solution of the following problem
\begin{equation*}
\begin{cases}
\vert\nabla u\vert=1 & \text{in $M\setminus K$}\,,\\
u=0 & \text{on $\partial K$}
\end{cases}
\end{equation*}
and we use the property of {\em semiconcavity} shared by such solutions to
obtain a rectifiability result for $Sing$.

Then, we investigate under which hypotheses also the closure of $Sing$ 
is rectifiable. This problem is strictly connected to the analysis of
the geodesic flow on $(M,g)$ originating from $K$, hence, we adapt 
ideas coming from the study of  the {\em cut locus} of a point in a
Riemannian manifold, which is actually the very special case when $K$
is a single point of $M$. 

Our results lead to the conclusion that, under some conditions on the
regularity of the set $K$, the Hausdorff dimension of the
closure of the singular set is at most  $(n-1)$ and that the gradient
of the distance function from $K$ is locally a vector field with {\em
  special bounded variation} (see~\cite{amb1,amb3,amb4}).\\
Moreover, we also study when the singular set shares an higher
regularity and we analyse in detail its topological structure if $M$
is a two--dimensional analytic surface and $K$ an analytic subset.

The study of the distance function and of the associate eikonal
equation $\vert\nabla u\vert=1$ is a special example of a connection
which can be extended to a large class of stationary
Hamilton--Jacobi equations. In the last section we discuss some
problems about the structure of the singular set of more general
viscosity solutions, suggested by some geometric results for the cut
locus of a point.

\section{Stationary Hamilton--Jacobi Equations on Manifolds}

Let $M$ be a smooth and connected, $n$--dimensional differentiable
manifold.

We consider the following Hamilton--Jacobi problem in $\Omega\subset M$,
\begin{equation*}
\begin{cases}
{\mathrm H}(x,du(x),u(x))=0&\text{in $\Omega$}\,,\\
u=u_0&\text{on $\partial\Omega$}
\end{cases}
\end{equation*}
where ${\mathrm H}:T^*\Omega\times\R\to\R$ and $T^*$ denotes the cotangent
bundle.

\def\testfun{\varphi}
\begin{dfnz}
Given a continuous function $u:\Omega \to \rl$ and a point $x \in M$, the
{\em superdifferential} of $u$ at $x$ is the subset of $T_x^*M$ defined by
$$
\supdiff u(x)=\left\{  d\testfun(x) \,|\,   \testfun\in C^1(M),
  \testfun(x)-u(x)=\min_M \testfun-u \right\}\,.
$$
\def\testfun{\psi}
Similarly, the set
$$
\subdiff u(x)=\left\{ d\testfun(x) \,|\, \testfun\in C^1(M),
  \testfun(x)-u(x)=\max_M \testfun-u 
\right\}
$$
is called the {\em subdifferential} of $u$ at $y$.\\
Notice that it is equivalent to replace the $\max$ ($\min$) on all $M$
with the maximum (minimum) in an open neighborhood of $x$ in $M$.
\end{dfnz}

\def\testfun{\varphi}
It is easy to see that $\supdiff u(x)$ and $\subdiff u(x)$ are 
both nonempty if and only if $u$ is
differentiable at $x\in M$. In this case we have
$$
\supdiff u(x)=\subdiff u(x)=\{ d u(x) \}\,.
$$
We list here without proof some of the standard properties of the sub and
superdifferentials which will be needed later.

\begin{prop}\label{supdiff and diffeo}
If $\chart:N\lra M$ is a map between the smooth manifolds $N$ and 
$M$ which is $C^1$ around $x\in N$, then
\begin{equation*}
\supdiff(u\composed\chart)(x)\supset\supdiff
u(\chart(x))\composed d\chart(x)\,=\,\{v\composed d\chart (x) \,|\, v\in
\supdiff u(\chart(x))\}\,.
\end{equation*}
If $\chart$ is a local diffeomorphism near $x$, the inclusion becomes an 
equality. An analogous statement holds for $\subdiff $. 
\end{prop}

\begin{prop}\label{supdiff and composition}
If $\theta:\real\lra\real$ is a $C^1$ function such that
 $\dot\theta(u(x))\geq 0$, then
$$
\supdiff(\theta\composed u)(x)\supset
 d\theta(u(x))\composed \supdiff u(x)\,=\,
\{d\theta(u(x))\composed v \,|\, v\in \supdiff u(x)\}\,,
$$
similarly for $\subdiff $. If $\dot\theta(u(x))>0$ then the inclusion
is an equality.
\end{prop}

For a locally Lipschitz  function $u$ on a Riemannian manifold
$(M,g)$, $\supdiff u(x)$ and $\subdiff u(x)$
are compact convex sets, almost everywhere coinciding with
the differential of the function $u$, by Rademacher's Theorem.\\
For a generic continuous function $u$ we prove in the next proposition
that  $\supdiff u(x)$ and $\subdiff u(x)$ are not empty in a dense
subset.

\begin{prop}\label{semidiffcontinue}
Let $u:\Omega\to\R$ be a continuous function on an open subset
$\Omega$ of $M$. Then the subdifferential $\subdiff u(x)$
(the superdifferential $\supdiff u(x)$) is not empty for every $x$ in a
dense subset of $\Omega$.
\end{prop}

\begin{proof}
It is always possible to endow $M$ with a Riemannian structure giving 
a metric $d(\cdot\,,\cdot)$ on $M$ which generates the same
topology.\\
Consider a generic point $y\in\Omega$ and a geodesic ball $B$
contained in $\Omega$ with center $y$. If the ball $B$ is small enough, the
function $x\mapsto d^2(x,y)$ is smooth in $\overline B$. Taking a large
positive constant $A$, the function $F_A(x)=u(x)+Ad^2(x,y)$ has a
local minimum at a point $x_A$ in the interior of $B$. At $x_A$ the
subdifferential of the function $F_A$ must contain the origin of
$T_{x_A}^*M$, hence, being $d^2(x,y)$  differentiable in the ball $B$,
the differential of $-d^2(x,y)$ at $x_A$ belongs to $\subdiff
u(x_A)$. As the point $y$ and the ball $B$ were arbitrarily chosen,
the set of points where the subdifferential of $u$ is not empty is 
dense in $\Omega$.\\
The same argument holds for the superdifferential of $u$, considering
the function $-u$.
\end{proof}

Now we introduce the notion of semiconcavity which will play a central
role in the first part of the paper.

\begin{dfnz}\label{semicsur}
Given an open set $\Omega\subset\rl^n$, a continuous function 
$u:\Omega\to\rl$ is called {\em locally semiconcave} if, for any open
convex set $\Omega^\prime\subset\Omega$ with compact closure in $\Omega$, there
exists a constant $C$ such that one of the following three equivalent
conditions is satisfied,
\begin{enumerate}
\item  $\forall x,h$ with $x$, $x+h$, $x-h\in \Omega^\prime$,
\begin{equation*}
u(x+h)+u(x-h)-2u(x) \leq 2C|h|^2\,,
\end{equation*}
\item  $u(x)-C\vert x\vert^2$ is a concave function in $\Omega^\prime$,
\item  $D^2u\le2C\,Id$ in $\Omega^\prime$, as distributions  ($Id$ is
  the $n\times n$ identity matrix).
\end{enumerate}
\end{dfnz}

In order to give a meaning to the concept of semiconcavity when the
ambient space is a differentiable manifold $M$, we analyse the
stability of this property under composition with $C^2$ maps.

\begin{prop} Let $\Omega$ and $\Omega^\prime$ two open subsets of
  $\R^n$. If $u:\Omega\to\rl$ is a Lipschitz function 
  such that $u(x)-C\,\vert x\vert^2$ is concave and
  $\chart:\Omega^\prime\to\Omega$ is a $C^2$ function with bounded first
  and second derivatives, then $u\composed\chart:\Omega^\prime\to\R$
  is a Lipschitz function and $u\composed\chart(y)-C^\prime|y|^2$ is 
  concave, for a suitable constant $C^\prime$.
\end{prop}

The proof is straightforward.\\
Then, the following definition is well--posed. 

\begin{dfnz}\label{semiconcave on M} A continuous function $u:M \to
\rl$ is called {\em locally semiconcave} if, for any local chart
$\chart:\real^n\lra \Omega\subset M$, the function $u\composed\chart$
is locally semiconcave in $\R^n$.
\end{dfnz}

The importance of semiconcave functions in connection with the 
generalized differentials is expressed by the following proposition
(see~\cite{PC-HMS:OnTheSin}).

\begin{prop}\label{supcontinu}
Let the function $u:M \to \rl$ be locally semiconcave, then the
superdifferential $\supdiff u$ is not empty at each point, moreover, 
$\supdiff v$ is {\em upper semicontinuous}, namely
$$
x_k\lra x ,\quad
  v_k\lra v,\quad  v_k\in \supdiff u(x_k)\quad \Longrightarrow \quad
  v\in \supdiff u(x)\,.
$$
In particular, if the differential $du$ exists at {\em every} point
of $\Omega\in M$, then $u\in C^1(\Omega)$.
\end{prop}

Now we introduce the definition of viscosity solution.\\
Let $\Omega$ be an open subset of $M$ and ${\mathrm H}$, called {\em
  Hamiltonian} function, a continuous real function on
$T^*\Omega\times\real$. We are interested in the following
Hamilton--Jacobi problem
\begin{equation}\label{HJ on M}
{\mathrm H}(x,du(x),u(x))=0\quad\quad\text{in $\Omega$}\,.
\end{equation}

\begin{dfnz}\label{viscodef} We say that a continuous function $u$ is
a {\em viscosity solution} of equation~\eqref{HJ on M} if for every 
$x\in\Omega$, 
\begin{equation}\label{visco HJ on M}
\begin{cases}
{\mathrm H}(x,v,u(x))\leq 0&\forall v\in \supdiff u(x)\,,\\
{\mathrm H}(x,v,u(x))\geq 0&\forall v\in \subdiff u(x)\,.
\end{cases}
\end{equation}
If only the first condition is satisfied (resp. the second) 
$u$ is called a {\em viscosity subsolution} (resp. a {\em viscosity
  supersolution}).
\end{dfnz}

If $\Omega^\prime$ is an open subset of another smooth differentiable
manifold $N$ and $\chart:\Omega^\prime\to\Omega$ is a $C^1$ local
diffeomorphism, we define the {\em pull--back} of the Hamiltonian
function $\chart^*{\mathrm H}:T^*\Omega^\prime\times\real\to\real$ by
\begin{equation*}
\chart^*{\mathrm H}(y,v,r)={\mathrm H}(\chart(y),v\composed d\chart(y)^{-1},r)\,.
\end{equation*}
Taking into account Proposition~\ref{supdiff and diffeo}, the
following statement is obvious.

\begin{prop}\label{visco and diffeo}
If $u$ is a viscosity solution of ${\mathrm H}=0$ in $\Omega\subset M$ and
$\psi:\Omega^\prime\to\Omega$ is a $C^1$ local diffeomorphism, then
$u\composed\chart$ is a viscosity solution of $\chart^*{\mathrm H}=0$ in
$\Omega^\prime\subset N$.
\end{prop}

\section{The Distance Function from a Subset of a Manifold}

From now on, $(M,g)$ will be a smooth, connected and complete
Riemannian manifold without boundary, of dimension $n$.

We will study the distance function $d_K$ from a closed and not empty
subset $K$ of  $M$ (for technical reasons, sometimes we consider also 
its square $d_K^2(x):M\to\R$).\\
The distance between two points $x$ and $y$ or from the point $x$ to
the set $K$ is defined as the infimum of the lengths of the $C^1$
curves starting at $x$ and ending at $y$, or on $K$, respectively. 
As $M$ is complete, by the Theorem of Hopf--Rinow, such infimum is
reached by at least one curve which will be a smooth geodesic.

The distance from the set $K$ is a continuous function on $M$ but in
general it is not everywhere differentiable, for instance, 
if the manifold $M$ is compact, the distance
function from any proper subset will be singular at the points of
absolute maximum. This section deals precisely with the set where the
gradient of $d_K$ does not exist.

In the following we will denote the distance between two points 
$x, y\in M$ with $d(x,y)$ and the {\em exponential map} of $(M,g)$ with
${\mathrm{Exp}}:TM\times\R\to M$. For simplicity, we will write 
$\vert v\vert$ for the modulus of a vector $v\in TM$, defined as
$\sqrt{g(v,v)}$.\\
Moreover, we will take often into account the identification between
the differential $du$ and the gradient $\nabla u$ of a function, 
induced by the scalar product $g$.

\begin{teo} The distance function $d_K$ is the {\em unique} 
viscosity solution of  the following Hamilton--Jacobi problem
\begin{equation}\label{H per dist}
  \begin{cases}
    \vert\nabla u\vert-1=0&\text{in $M\setminus K$}\,,\\
    u=0&\text{on $K$}
  \end{cases}
\end{equation}
in the class of continuous functions
bounded from below.\\
The function $d_K^2/2$ is the {\em unique} viscosity solution of
  \begin{equation}\label{distprob}
    \begin{cases}
      \vert\nabla u\vert-2u=0&\text{in $M$}\,,\\
      u=0&\text{on $K$}
    \end{cases}
  \end{equation}
in the class of continuous functions on $M$
such that their zero set is $K$.
\end{teo}

\begin{rem} The restriction to lower bounded functions is necessary,
  $\Vert x\Vert$ and $-\Vert x\Vert$ are both viscosity solutions of
  Problem~\eqref{H per dist} with $M=\R^n$ and $K=\{0\}$. Moreover, the
  completeness of $M$ plays an important role here, if $M$ is the open
  unit ball of $\R^n$ the same example shows that the uniqueness does not
  hold.\\
  Notice also that every function $d_H^2/2$ where $H$ is a
  closed subset of $M$ with $H\supset K$, is a viscosity solution of
  Problem~\eqref{distprob}, equal to zero on $K$.
\end{rem}

\begin{proof}
\def\uzero{u_0}
\def\e{\varepsilon}
\def\y{{\widetilde x}}
\def\vv{d_K}
Notice that  $d_K(x)$ is the minimum time $t\ge 0$
 for any curve $\gamma$ to reach
a point $\gamma(t)\in K$, subject to the conditions
 $\gamma(0)=0$ and $|\dot\gamma|\le 1$;
$d_K$ is then the value function of a ``minimum time problem'';
this proves that $d_K$ is also a viscosity solution of 
Problem~\eqref{H per dist}, by well known results (see for example
Proposition~2.3, Chapter~IV in~\cite{b-cd}). Then 
we show that the function $d_K^2/2$ is a solution of 
Problem~\eqref{distprob}.\\ 
First of all, notice that the distance function from $K$ is a
1--Lipschitz function, hence $d_K^2$ is locally Lipschitz.\\
As $d_K$ is 1--Lipschitz, at every point of $K$ the
function $d^2_K$ is differentiable and its differential is 
zero. Hence, the definition of viscosity solution holds for 
points belonging to $K$. In order to prove the theorem, it is then
sufficient to test conditions~\eqref{visco HJ on M} on the generalized
differentials at the points of the open set $M\setminus K$.\\
Since $d_K^2/2$ is positive in $M\setminus K$, applying
Proposition~\ref{supdiff and composition}  with the function
$\theta(t)=\sqrt{2t}$, we see that the function $d_K^2/2$ is a
viscosity solution of
\begin{equation*}
g\left(\frac{\nabla u}{\sqrt{2u}},\frac{\nabla u}{\sqrt{2u}}\right)-1=0
\end{equation*} 
in $M\setminus K$. Being there positive, it  also solves
\begin{equation*}
g(\nabla u,\nabla u)-2u=0
\end{equation*}
in $M\setminus K$. This fact together with the previous remark about 
the behavior of $d_K^2$ at the points of $K$ gives the claim.

\def\p{{\Psi}}
\def\Bc{\delta \,}
\def\e{\varepsilon}
Suppose now that $u$ is a viscosity
solution of Problem~\eqref{H per dist} then, $u$ is also a solution of
\begin{equation*}
  \begin{cases}
    \vert\nabla u\vert-1=0&\text{in $M\setminus K$}\,,\\
    u=0&\text{on $K$}\,.
  \end{cases}
\end{equation*}
As in the work of
  Kru\v{z}hkov~\cite{kruz1}, we consider the function $v= -e^{-u}$
  which, by Proposition~\ref{supdiff and composition}, turns out to be
  a viscosity solution of
\begin{equation*}\label{nuova}
\begin{cases}
|\nabla v|+v =0 & \text{ in $M\setminus\bordreg$}\,,\\
 v=-1& \text{ on $\bordreg$}
\end{cases}
\end{equation*}
moreover, $|v| \le e^{-\inf u}$.\\
We establish an uniqueness result for this last problem in the class
of bounded functions $v$, which clearly implies the
first uniqueness result. We remark that the proof is based on
similar ones in~\cite{crevli1,crisli1,YG-SG-HI-MHS:ComPriAnd}.\\
We argue by contradiction, suppose that $u$ and $v$ are two bounded
solutions of~\eqref{nuova}, $\vert u\vert$, $\vert v\vert\leq C$, 
and that at a point $\overbar x$ we have $u(\overbar x)\geq
2\e+v(\overbar x)$ with $\e >0$.\\
Let $b(x,y):M\times M \lra \real$ be a smooth function satisfying
\begin{itemize} 

\item $b\ge 0$ 

\item $\left\vert\nabla_x b (x,y) \right\vert$, $\left\vert\nabla_y b (x,y) \right\vert \le 2$

\item $\sup_{M\times M} \left\vert d(x,y)-b(x,y)\right\vert < \infinity$
\end{itemize}
such a function can be obtained {\em smoothing} the distance function
in $M\times M$.\\
We fix a point $x_0$ in $K$ and we define the smooth function $B(x) =
b(x,x_0)^2$. By the properties of $b$ and the boundedness of $u$ and
$v$, the following function $\Psi:M\times M\to\R$
\begin{equation*}
\p(x,y) = u(x) - v(y) - \l d(x,y)^2-\Bc B(x) - \Bc B(y)
\end{equation*}
has a maximum at a point $\widehat x,\widehat y$ (dependent on
the positive parameters $\delta$ and $\l$) and such maximum $\p(\widehat
x,\widehat y)$ is less than $2C$. Hence, the function
\begin{equation}\label{funz1}
x\mapsto [ v(\widehat y) + \l{d(x,\widehat y)^2} + \Bc B(x) +\Bc B(\widehat y)] -
u(x)
\end{equation}
has a minimum at $\widehat x$ while
\begin{equation}\label{funz2}
y\mapsto [ u(\widehat x) - \l{d(\widehat x,y)^2} - \Bc B(\widehat x) -\Bc B(y)] -
v(y)
\end{equation}
has a maximum at $\widehat y$.\\
If $2\delta \le \e / B(\overbar x)$ then
$$
\p(\widehat x,\widehat y)\ge\p(\overbar x,\overbar x) \ge
2\e-2\delta B(\overbar x)\ge\e
$$
hence, we get 
\begin{equation}\label{funz3}
\delta B(\widehat x)+ \delta B (\widehat y)+ \l d(\widehat x,\widehat y)^2+ \e \le 
u(\widehat x) -v(\widehat y) \le 2C\,.
\end{equation}
This shows that, for a fixed $\delta$, the pair ${\widehat
  x},{\widehat y}$ is contained in a bounded set and, if $\lambda$
  goes to $+\infty$ the distance between $\widehat x$ and $\widehat y$
  goes to zero. Possibly passing to a subsequence for $\lambda$
  going to infinity, ${\widehat x}$ and ${\widehat y}$ converge to a
  common limit point $z$ which cannot belong to $K$, otherwise we
  would get $\varepsilon\leq u(z)-v(z)=0$, thus, for some $\lambda$
  large enough also ${\widehat x}$ and ${\widehat y}$ do not belong to
  $K$.\\
  As the function $d^2(x,y)$ is smooth in $B_z\times B_z\subset
  M\times M$, where $B_z$ is a small geodesic ball around $z$,
  choosing a suitable $\lambda$ large enough we can differentiate the
  functions inside the square brackets in equations~\eqref{funz1}
  and~\eqref{funz2} obtaining
  $$
  \widehat v=\Bc\nabla B(\widehat x) + \l\nabla_x d^2 (\widehat x,\widehat y) ~~~~\in
  \supdiff  u(\widehat x)\,,
  $$
$$\widehat w=- \Bc \nabla B(\widehat y) - \l\nabla_y d^2 (\widehat x,\widehat y) ~~~\in
\subdiff v(\widehat y)\,.
$$
By Definition~\ref{viscodef} we have that $ |\widehat v| +u(\widehat x)\le 0$ and
$ |\widehat w| + v(\widehat y)\ge 0$, hence
$$
u(\widehat x) -v(\widehat y) + |\widehat v| - |\widehat w| \le 0\,.
$$
Moreover,
\begin{align*}
  |\widehat v|-|\widehat w| =&\left\vert  \delta\, \nabla B(\widehat v) + \l \nabla_x d^2
    (\widehat x,\widehat y)\right\vert -\left\vert\delta\, 
    \nabla B(\widehat y) + \l\nabla_y d^2
    (\widehat x,\widehat y)\right\vert\\
  \ge&\left\vert\l\nabla_x d^2 (\widehat x,\widehat
    y)\right\vert-\left\vert\l\nabla_y d^2 (\widehat x,\widehat y)\right\vert -
  |\delta\, \nabla B(\widehat y)| - |  \delta\, \nabla B(\widehat x)  |\\
=&\,2\lambda d(\widehat x,\widehat y)\left\vert\nabla_x d(\widehat
  x,\widehat y)\right\vert -2\lambda d(\widehat x,\widehat y)\left\vert\nabla_y d(\widehat
  x,\widehat y)\right\vert-|\delta\, \nabla B(\widehat y)| - |
\delta\, \nabla B(\widehat x)  |\\
=&\,2\lambda d(\widehat x,\widehat y)-2\lambda d(\widehat x,\widehat y)-|\delta\, \nabla B(\widehat y)| - |
\delta\, \nabla B(\widehat x)  |\\
=&-|\delta\, \nabla B(\widehat y)| - |\delta\, \nabla B(\widehat x)
|
\end{align*}
which implies,
\begin{equation*}
    u(\widehat x) -v(\widehat y)-\delta |\nabla B(\widehat y)| -
    \delta| \nabla B(\widehat x)|\leq0\,.
\end{equation*}
Finally, we have that 
$$
\delta |\nabla B(\widehat x)|=2\delta|b(\widehat x,x_0)\nabla
b(\widehat x,x_0)|\leq4\delta\sqrt{B(\widehat x)}
$$
and using the estimate $\delta B(\widehat x)\leq 2C$ which follows
from equation~\eqref{funz3},
$$
\delta |\nabla B(\widehat x)|\leq8\sqrt{2\delta C}\leq\e/4
$$
if $\delta$ was chosen small enough. Holding the same for $\widehat
y$, we conclude that
$$
u(\widehat x) - v( \widehat y ) - \e/2\leq0
$$
which is in contradiction with the fact that $u(\widehat x)
-v(\widehat y)\ge \e$.

About the second uniqueness claim, if $u$ is a 
continuous viscosity solution of Problem~\eqref{distprob} then, by
Proposition~\ref{semidiffcontinue} the superdifferential of $u$ is not 
empty in a dense subset of $M\setminus K$, hence, directly by
the equation and by continuity, $u$ is non negative. By the hypothesis
on its zero set we conclude that $u$ is positive in all $M\setminus
K$. Composing $u$ with the function $t\mapsto \sqrt{2t}$, we see that
$\sqrt{2u}$ is a positive, continuous viscosity solution of
Problem~\eqref{H per dist}, then it must coincide with $d_K$, by the
previous result. It follows that $u=d^2_K/2$.
\end{proof}

We now study the rectifiability of the {\em singular set} of $d_K$, 
\begin{equation*}
Sing=\left\{x\in M\,\vert\,\text{$d_K^2$ is not differentiable at
  $x$}\right\}\,.
\end{equation*}

\begin{rem}
In this definition we used the squared distance function instead of 
the distance in order to avoid to consider also the points of
the boundary of $K$, which are singular for $d_K$. It is easy to see
that outside $K$ the distance and its square have the same regularity.
\end{rem}

\begin{prop}\label{carlosemicref}
The function $d_K$ is locally semiconcave in $M\setminus K$.
\end{prop}

\begin{proof}
  The distance function $d_K$ is a viscosity solution of ${\mathrm H}=0$ in
  $M\setminus K$, where the Hamiltonian function is given by
  ${\mathrm H}(x,v,t)=|v|^2-1$. We choose a smooth local chart
  $\chart:\real^n\to\Omega\subset M$ and we define
  $v=d_K\composed\chart$, which is a locally Lipschitz function and,
  by Proposition~\ref{visco and diffeo}, it is a viscosity solution of 
  $\chart^*{\mathrm H}=0$.\\
  The pull--back of the Hamiltonian function on $\R^n$ takes the form
  \begin{equation*}
    \chart^*{\mathrm H}(y,w,s)=g_{\chart(y)}(d\chart(w),d\chart(w))-1=g_{ij}(y)w_iw_j-1
  \end{equation*}
  for $(y,w,s)\in\R^n\times\R^n\times\R$ and where $g_{ij}(y)$ are the
  components of the metric tensor of $M$ in local coordinates.\\
  Since the matrix $g_{ij}(y)$ is positive definite $\chart^*{\mathrm H}(y,w,s)$
  is locally uniformly convex in $w$, hence, by Theorem~5.3
  of~\cite{PLL:GenSolOf}, it follows that $v=d_K\composed\chart$ is
  locally semiconcave in $\R^n$. Recalling Definition~\ref{semiconcave on M},
  this means that $d_K$ is locally semiconcave in $M\setminus K$.
\end{proof}

The semiconcavity of $d_K$ allows us to work with the
superdifferentials when the gradient does not exist. Indeed, notice 
that the points of $Sing$ are precisely those where the
superdifferential is not a singleton, then, the following result is a
straightforward consequence of Proposition~\ref{supcontinu}.

\begin{prop}\label{c1reg} The function $d^2_K$ is of class $C^1$ in
  the open set $M\setminus\overline{Sing}$ and $d_K$ is $C^1$ in
  $M\setminus\left(K\cup\overline{Sing}\right)$.
\end{prop}

\begin{rem}\label{estremali} The semiconcavity property also gives
  information about the  propagation of singularities and the
  relations between the structure of the superdifferential at a point
  $x$ and the set of minimal geodesics from $x$ to $K$ (see
  ~\cite{alambca,ambcaso}).\\
The set $Ext(\partial^+d_K^2(x)/2)$ of
extremal points of the (convex) superdifferential set of 
$d_K^2/2$ at $x$ is in one--to--one correspondence with the family
${\mathcal G}(x)$ of minimal geodesics from $x$ to $K$. Precisely
${\mathcal G}(x)$ is described by
\begin{equation}\label{geomin}
  {\mathcal G}(x)=\left\{{\mathrm{Exp}}(x,-v,\cdot):[0,1]\lra M
    \,|\,\text{for $v\in
      Ext(\partial^+d_K^2(x)/2)$}\right\}\,.
\end{equation}
The set of points of $K$ at minimum distance from $x$ are given by
${\mathrm{Exp}}(x,-v,1)$ for $v$ in the set of extremal points of the
superdifferential set of $d_K^2/2$ at $x$. 
As a particular case we have that if the function $d_K^2$ is differentiable
  at $x$, then the point of $K$ closest to $x$ is uniquely
  determined and given by ${\mathrm {Exp}}(x,-\nabla d_K^2(x)/2,1)$.\\
Finally, notice that $Sing$ is the set of points $x$ such that 
the distance $d(x,K)$ is realized by more than one minimal geodesic
between $x$ and $K$.
\end{rem}

The rectifiability of $Sing$ now follows. 

\begin{prop}\label{rect1}
The set $Sing$ is $C^2$--rectifiable.
\end{prop}

\begin{proof}
By a result proved in~\cite{alberti1}, the singular set of a
locally semiconcave function in an open set of $\R^n$ is
$C^2$--rectifiable. We take a countable family of local charts
$\chart_i:\R^n\to\Omega_i$ and consider the functions
$d_K\composed\chart_i$. These functions are locally semiconcave in
$\R^n$ with singular sets $Sing_i$, hence, by the relation
$$
Sing\subset\bigcup_{i=1}^{\infty}\chart_i(Sing_i)
$$
we get the thesis.
\end{proof}

The same statement does not hold for the closure of $Sing$, for a
generic closed set. We describe now a counterexample
showing that indeed the set $\overbar{Sing}$ is not rectifiable 
for a set $K$ of class $C^{1,1}$ only.

We look for a convex open set $\Omega$ with a $C^{1,1}$
boundary in $\real^2$ such that the closure of the singular set
$Sing$ of the distance function from its boundary has nonzero Lebesgue
measure, hence it is not rectifiable.\\
We start with a Cantor--like set ${\mathcal C}\subset\SS^1$, closed 
with empty interior in $\SS^1$, with no isolated points and positive
Hausdorff ${\mathcal H}^1$ measure.
Such a set can be constructed as follows
\begin{equation*}
{\mathcal C}=\SS^1\setminus\bigcup_{i=1}^\infty I_i
\end{equation*}
where $\{I_i\}$ is a countable family of open disjoint connected arcs
on $\SS^1$, whose middle points are $p_i\in\SS^1$ and such that the sum
of their lengths is less than $2\pi$.\\
We claim that every point of ${\mathcal C}$ is a limit point of the sequence
$\{p_i\}$. If $p\in{\mathcal C}$ there must be a sequence
of arcs $I_{i_j}$ arbitrarily close to $p$, since the arcs are countable
and the sum of their lengths is bounded by $2\pi$ we have that they
shrink when $j$ goes to infinity, hence $p_{i_j}\to p$.\\
We define an open convex set ${\Omega^\prime}$ as the intersection of
the open halfplanes, containing the origin of $\R^2$, determined by
the tangent lines to $\SS^1$ at the points of ${\mathcal C}$, see the
following figure.\\
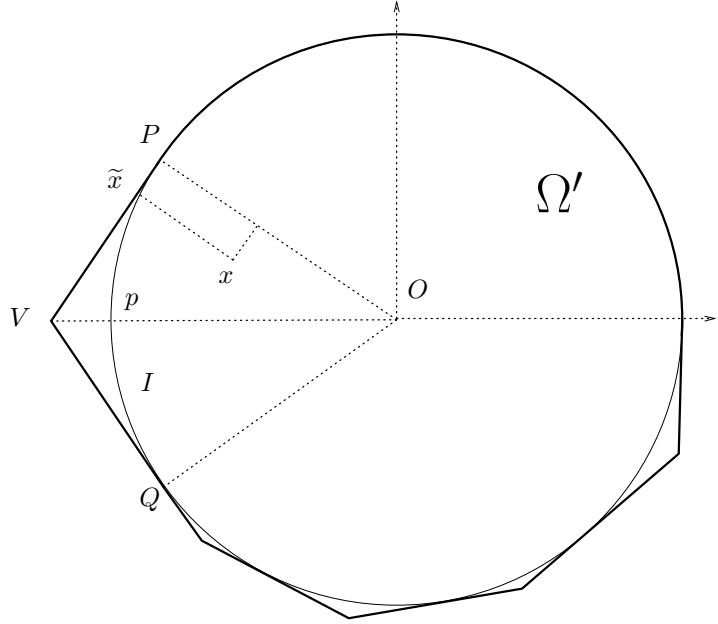
\begin{figure}[ht]
\setlength{\unitlength}{0.00625in}
\begin{picture}(590,531)(0,-10)
\thicklines
\put(324,250){\arc{478}{3.7144}{6.3388}}
\thinlines\thinlines
\put(324,250){\ellipse{478}{478}}
\dottedline{5}(324,251)(324,516)
\path(326.000,508.000)(324.000,516.000)(322.000,508.000)
\dottedline{5}(324,251)(590,251)
\path(582.000,249.000)(590.000,251.000)(582.000,253.000)
\dottedline{5}(207,328)(187,300)(110,354)
\thicklines
\path(125,382)(35,249)(129,110)
        (161,65)(284,0)(429,25)
        (560,138)(563,250)
\thinlines
\dottedline{5}(324,250)(129,110)
\dottedline{5}(324,250)(129,382)
\dottedline{5}(324,250)(35,249)
\put(110,190){\makebox(0,0)[lb]{\smash{{{$I$}}}}}
\put(97,262){\makebox(0,0)[lb]{\smash{{{$p$}}}}}
\put(109,398){\makebox(0,0)[lb]{\smash{{{$P$}}}}}
\put(0,245){\makebox(0,0)[lb]{\smash{{{$V$}}}}}
\put(440,340){\makebox(0,0)[lb]{\smash{{{\Huge{$\Omega^\prime$}}}}}}
\put(333,268){\makebox(0,0)[lb]{\smash{{{$O$}}}}}
\put(82,360){\makebox(0,0)[lb]{\smash{{{$\widetilde{x}$}}}}}
\put(109,94){\makebox(0,0)[lb]{\smash{{{$Q$}}}}}
\put(175,281){\makebox(0,0)[lb]{\smash{{{$x$}}}}}
\end{picture}
\end{figure}\\
Let us take an arc $I$ with middle point $p$, bounded by
$P,\,Q\in{\mathcal C}$ and consider the associate quadrilateral
$OPVQ$. If the point $x$ is inside the
open triangle $OPV$ it is clear that the point of $\partial
\Omega^\prime$ closest to $x$ belongs to the segment ${PV}$ and it is
unique ($\widetilde{x}$ in the figure). Hence, for such points the
distance from the boundary of $\Omega^\prime$ coincide with the distance from
the segment ${PV}$.\\
Applying the same argument to the open triangle $OQV$, we see that the
segment ${OV}$ consists of singular points of $d_{\partial\Omega^\prime}$,
moreover, the segment ${OV}$ intersects $\SS^1$ at the point $p$.\\
It follows that the union ${\mathcal S}$ of the segments from the
middle points $p_i$ to the origin coincides with $Sing$ for
$d_{\partial\Omega^\prime}(x)$. Being $p_i$ dense in ${\mathcal C}$,
the closure of ${\mathcal S}$ contains $\lambda{\mathcal
  C}\subset\R^2$ for every $\lambda\in[0,1]$. As ${\mathcal C}$ has
$\H^1$ positive measure, the Lebesgue measure of ${\overbar{\mathcal
    S}}$ is positive, by Fubini's Theorem.\\
Now let $\Omega$ be the set of points of $\R^2$ whose distance from
the convex $\Omega^\prime$ is less than 1. It is immediate to check
that
$$
d_{\partial\Omega}(x)=d_{\partial\Omega^\prime}(x)+1\quad\quad\forall
  x\in{\Omega^\prime}
$$
hence for every $x$ in the unit ball.\\
So the closure of $Sing$ for the distance function from the
boundary of $\Omega$ (or from the complementary set of $\Omega$ in
$\R^2$) has positive Lebesgue measure, moreover, by the properties of
convex bodies, the boundary of $\Omega$ is of class at least
$C^{1,1}$. Since the Lebesgue measure of $\overline{Sing}$ is positive
it cannot be rectifiable.

\begin{rem}
In the next section we will show that if the boundary of $K$ is of class
at least $C^3$ then also the closure of $Sing$ is rectifiable. To our
knowledge it is unknown even in $\R^2$ if the gap between such result
and the previous counterexample can be filled, that is, if the $C^2$
(or maybe $C^{2,1}$) regularity of the boundary of $K$ is enough to
get the rectifiability of the closure of the singular set.
\end{rem}

\section{Rectifiability of the Closure of the Singular Set}

In this section we are going to show that an higher regularity of the
set $K$ implies the rectifiability also of the closure of the singular
set. Moreover, we determine a relation between the regularity of $K$
and of the hypersurfaces covering $Sing$.

In all this section $K$ is a $k$--dimensional embedded
$C^r$ submanifold of $M$ without boundary, with $0\leq k\leq n-1$ (the
case $k=n$ is trivial) and $r\ge 2$.\\
Let $UK$ be the unit normal bundle of $K$ in $M$, we denote with 
$F:UK\times[0,+\infty)\to M$ the restriction of the exponential map of
$M$ to $UK\times[0,+\infty)$. Since $K$ is $C^r$, $UK$ is a manifold
of class $C^{r-1}$ and being ${\mathrm {Exp}}$ a smooth map,
$F(q,v,t)$ and all its derivatives in the $t$ variable are $C^{r-1}$
functions.

\begin{rem}
The case when $K$ is the closure of an open set of $M$ with smooth
boundary can be reduced to our case. Indeed, 
if $x\in M \setminus K$, then the minimal geodesic
from $x$ to $K$ ends in $\partial K$ without touching the interior
points, hence to study $d_K$ we can simply consider the
distance function from $\partial K$ in every open connected component
of $M\setminus K$. However, in this case, some results like the higher
smoothness of $d_K^2$ at the points of $\partial K$ expressed by
Proposition~\ref{regular}, could be lost.
\end{rem}

The behavior of $d_K$ near $K$ is well known.

\begin{prop}\label{regular}
There exist $\varepsilon>0$ and an open neighborhood
  $\Omega$ of $K$ in $M$ such that the map 
  $F\vert_{UK\times(0,\varepsilon)}:{UK\times(0,\varepsilon)}\to
  \Omega\setminus K$ is a $C^{r-1}$ diffeomorphism.\\
Moreover,
  \begin{itemize}
  \item for every point in $\Omega$ there is an unique point of
    minimum distance in $K$ (the {\em unique projection property}
    holds for $K$ in $\Omega$),
\item the distance function $d_K$ is $C^r$ in $\Omega\setminus K$, 
  \item the squared distance function $d^2_K$ is $C^r$ in $\Omega$.
  \end{itemize}
\end{prop}

\begin{rem} It can be proved that $K$ has to be at least $C^{1,1}$ in
  order to share the unique projection property in a neighborhood, 
  in such case the squared distance function also turns out to be of
  class $C^{1,1}$.\\
  See~\cite{delzol1,delzol2} for a detailed discussion of the
  relation between the regularity of $K$ and of $d_K$.
\end{rem}

In order to study what happens {\em far} from $K$ 
we have to analyse the sets of points where the {\em unique
  projection property} fails or $F$ is not a local diffeomorphism. 
From a topological point of view, the problem is
naturally connected with the study of the singularities of maps
between Euclidean spaces. For instance, when $K$ is a single
point of $M$ the singular sets were shown to be related to the
classes of singularities considered by the {\em Theory of Catastrophes},
see~\cite{buchner}.

Consider the geodesic curve $t\mapsto F(q,v,t)$ for $t\in [0,t_0]$
($(q,v)\in UK$ is fixed), for small values of $t_0$ it is the unique
minimizer of the length functional between its end point
$p=F(q,v,t_0)$ and $K$ but for large $t_0>0$ it could cease to be
minimal. Hence, there exists a value $\sigma$ (possibly $+\infty$)
such that this geodesic is minimal between $q$ and $F(q,v,t)$ for
every $t\leq\sigma$, but not on any larger interval. 
If $\sigma(q,v)<+\infty$, we say that the point $F(q,v,\sigma(q,v))$
is the {\em cut point} of the geodesic $F(q,v,t)$ and we define the
following set,
\begin{equation*}
V_K=\left\{(q,v,t)\in
UK\times\R^+\,\vert\,\text{$t<\sigma(q,v)$}\right\}\,.
\end{equation*}
Notice that the set $F(V_K)$ clearly contains $\Omega\setminus K$,
where $\Omega$ is the open neighborhood of Proposition~\ref{regular}.

\begin{dfnz}
The set of points $F(q,v,\sigma(q,v))$ for $(q,v)\in UK$ with 
$\sigma(q,v)<+\infty$ is called the {\em cut locus} of $K$, we denote
it with $\cutt$.
\end{dfnz}

The reasons why a geodesic ceases to be minimal are explained by the
following proposition (see~\cite{gahula,sakai}).

\begin{prop}\label{carlo100}
If for a geodesic $F(q,v,t)$ we have $\sigma(q,v)<+\infty$, at least one of the 
following two non exclusive conditions is satisfied:
\begin{enumerate}
\item at the point $p=F(q,v,\sigma(q,v))$ there arrives another minimal 
  geodesic from $K$,
\item the differential $dF(q,v,\sigma(q,v))$ is not invertible.
\end{enumerate}
Conversely, if at least one of these conditions is satisfied the 
geodesic $F(q,v,t)$ cannot be minimal on an interval larger that
$[0,\sigma(q,v)]$.
\end{prop}

Notice that, by Remark~\ref{estremali}, if condition~1 above is
satisfied then the point $p$ belongs to $Sing$, while $Sing$ is clearly
included in $\cutt$.\\

Then, we consider the following two subsets of $\cutt$,
\begin{itemize}
\item $Sing$ is the set of points $p=F(q,v,\sigma(q,v))$ where more
  than one minimal geodesic from $K$ arrives,
\item $J$ is the set of points $p=F(q,v,\sigma(q,v))$ such that the
  differential $dF(q,v,\sigma(q,v))$ is not invertible.\\
We call $J$ the {\em locus of optimal focal points}.
\end{itemize}
Clearly $Sing\cup J=\cutt$.

The next proposition (see~\cite{gahula,sakai}) establishes 
the connection between $\cutt$ and the distance function from $K$.

\begin{prop}\label{chiusura} The following statements hold,
\begin{enumerate}
\item $\cutt=\overline{Sing}$, that is, the cut locus of $K$ is closed
in $M$ and $Sing$ is a dense subset.
\item The set $V_K$ is open in $UK\times\R^+$.
\item The map $\sigma:UK\to\R^+$ is continuous.
\item The map $F$ is a $C^{r-1}$ bijection between $V_K$ and
  $M\setminus(K\cup\cutt)$ with a $C^{r-1}$ inverse.
\item The cut locus $\cutt$ is equal to $F(\partial V_K)$ where the
  boundary is considered in the ambient space $UK\times\R^+$.
\item The set $M\setminus\cutt$ can be continuously retracted on $K$,
  and, if $\sigma(q,v)<+\infinity$ for every $(q,v)\in UK$, the
  set $M\setminus K$ can be continuously retracted on the cut locus
  $\cutt$.
\item The open set $M\setminus\cutt$ has the unique projection
  property, moreover the squared distance function $d_K^2$ is of class
  $C^r$ in it. The distance $d_K$ is $C^r$ in  $M\setminus(\cutt\cup
  K)$.
\end{enumerate}
\end{prop}

By the point~1, the closure of $Sing$ is precisely the cut locus of $K$,
which is the union of $Sing$ and $J$. We study then the
rectifiability of these two sets separately, as the rectifiability
property is clearly stable under countable union (see
Definition~\ref{rectdef}).

We have seen in Proposition~\ref{rect1} that the set $Sing$ is always 
$C^2$--rectifiable. We partially improve this result when $K$ is more
regular.

\begin{prop}\label{crregula} If $K$ is of class $C^r$ with $r\geq3$ the
  set $Sing\setminus J$ is a $C^r$--rectifiable subset of $M$.
\end{prop}

We need a preliminary lemma.

\begin{lemma}\label{infinite then focal}
If there are infinitely many minimal geodesics from $p$ to $K$, then
$p$ is an optimal focal point.
\end{lemma}

\begin{proof}
If $F(q_i,v_i,\sigma(q_i,v_i))=p$ for infinite distinct geodesics
$F(q_i,v_i,t)$, then $\sigma(q_i,v_i)=d(q_i,p)=d_K(p)$ hence, by compactness, we
may assume that $(q_i,v_i)\lra(q,v)$ for some $(q,v)\in UK$. It
follows (by the semicontinuity of the length functional) that
$F(q,v,t)$ is a minimal geodesic for $p$ and that $dF(q,v,t)$ is
singular since $F$ is not locally injective near $(q,v,t)$.
\end{proof}

\begin{proof}[Proof of Proposition~\ref{crregula}]
Let $p$ be a point in $Sing\setminus J$.
We know that the number of minimal geodesics $F(q_i,v_i,t)$ from
$p=F(q_i,v_i,\sigma(q_i,v_i))$ to $K$ is finite by the lemma above
and greater than one, by the singularity at $p$. Moreover, the
differential $dF(q_i,v_i,\sigma(q_i,v_i))$ is invertible for every
$i$, then $F$ is locally invertible in the neighborhood of every point 
$(q_i,v_i,\sigma(q_i,v_i))\in UK\times\R^+$.\\
Let $U_i$ be disjoint open neighborhoods of
$(q_i,v_i,\sigma(q_i,v_i))$ such that $F$ restricted to every $U_i$ is a 
$C^{r-1}$ diffeomorphism with its image. We can also
suppose that $F(U_i)=U$ where $U$ is an open neighborhood of $p$ in
$M$. We define the functions $d_i:U\to\R^+$ given by
\begin{equation*}
d_i(x)=\pi_{\R^+}\left(F^{-1}(x)\cap U_i\right)
\end{equation*}
where $\pi_{\R^+}:UK\times\R^+\lra\R^+$ denotes the projection on the
second factor. Applying Proposition~\ref{regular} to a small
neighborhood of every $q_i$, we can see that all the functions $d_i$
are of class $C^r$.

The singular set $Sing\cap U$ is clearly contained in the union of
the sets $S_{ij}=\left\{r\in U\,\vert\, d_i(x)=d_j(x)\right\}$ for
$i\not=j$. We now prove that such sets are locally $C^r$
hypersurfaces. By the implicit function theorem, it is sufficient to
show that $\nabla d_i(x)-\nabla d_j(x)\not=0$ at the points of $S_{ij}$.
If $r=F(s_i,w_i,t)=F(s_j,w_j,t)$, for $(s_i,w_i,t)\in U_i$ and
$(s_j,w_j,t)\in U_j$ then
\begin{equation*}
\nabla d_i(x)=\nabla d_j(x)\quad\Longrightarrow
\quad\frac{dF}{dt}(s_i,w_i,t)=\frac{dF}{dt}(s_j,w_j,t)\,.
\end{equation*}
Hence, by the uniqueness of the geodesic from $r$ with a certain
initial velocity vector, we would get $(s_i,w_i)=(s_j,w_j)$,
contradicting the hypothesis that $U_i$ and $U_j$ are disjoint.
\end{proof}

The next proposition gives the rectifiability of the set $J$ of
optimal focal points which implies the rectifiability of the cut locus
of $K$, by the previous discussion. 

\begin{prop}\label{conrect} If $K$ is of class $C^r$ with $r\geq3$,
  then the set $J$ is $C^{r-2}$--rectifiable.
\end{prop}

\begin{teo} If $K$ is of class $C^r$ with $r\geq3$,
the closure of $Sing$, that is, the cut locus of $K$, is
$C^{r-2}$--rectifiable.
\end{teo}

\begin{proof}[Proof of Proposition~\ref{conrect}] 
We introduce the set ${\widetilde J}$ of the {\em first focal
  points} as follows. 
Let $F(q,v,t)$ be a geodesic from $K$ with $(q,v)\in UK$ and
$t\in\R^+$, considering the first value $t=c(q,v)$ such that
$dF(q,v,t)$ is not invertible or setting $c(q,v)=+\infty$ if $dF(q,v,t)$
is invertible for every $t\in\R^+$, we define the map
$c:UK\to\R^+\cup\{+\infty\}$. If $c(q,v)<+\infty$ we say that
$F(q,v,c(q,v))$ is the first focal point along the geodesic $F(q,v,t)$.\\
We consider the following set of points $G$ in $UK\times\R^+$ where
$dF$ is not invertible
\begin{equation*}
G=\left\{(q,v,c(q,v))\in UK\times\R^+\,\vert\,
\text{$c(q,v)<+\infty$}\right\}
\end{equation*}
and we call ${\widetilde J}=F(G)$ {\em locus of the first focal
  points} of $K$.\\
By Proposition~\ref{carlo100}, we have that
${\widetilde J}\supset J$, the set of optimal focal points, hence, it 
sufficient to show that the set ${\widetilde J}$ is
rectifiable to prove the same for $J$ and conclude the proof.

At the points of the set $G$ the rank of $dF$ is at
most $(n-1)$. We split $G$ in two subsets,
\begin{equation*}
\begin{split}
G_1=\left\{(q,v,c(q,v))\in G\,\vert\,\text{ Rank~$dF(q,v,c(q,v))=n-1$}\right\}\\
G_2=\left\{(q,v,c(q,v))\in G\,\vert\,\text{ Rank~$dF(q,v,c(q,v))<n-1$}\right\}\,.
\end{split}
\end{equation*}

The following version of Sard's Theorem can be found
in the book of Federer~\cite[Theorem~3.4.3]{fede}.

\begin{lemma}
Let $F:\R^n \to \R^n$ be a map of class $C^l$ for some $l\geq 1$.\\
If we set for any $k\in\{0,1,\dots,n-1\}$
\begin{equation*}
A_k=\left\{x\in\R^n\,\vert\,\text{ {\rm Rank}~$dF(x)\leq
    k$}\right\}\,,
\end{equation*}
then the Hausdorff measure ${\mathcal H}^{k+\frac{n-k}{l}}$ of
$F(A_k)$ is zero.
\end{lemma}
Considering local charts of $UK\times\R^+$ and $M$ and applying this
lemma to our map $F$ which is of class $C^{r-1}$ we get that
${\mathcal H}^{n-2+2/(r-1)}(F(G_2))=0$, which implies
${\mathcal H}^{n-1}(F(G_2))=0$ as $r\geq 3$.

Now we show that the set $G_1$ is locally a
$C^{r-2}$ hypersurface in $UK\times\R^+$ using the implicit function
theorem, this proves that $F(G_1)$ is a $C^{r-2}$--rectifiable
set in $M$ and so ${\widetilde J}$.

Let $(q,v,c(q,v))$ be a point in $G_1$, by the lower semicontinuity of
the rank, choosing a small neighborhood $B$ of $(q,v,c(q,v))$ in
$UK\times\R^+$ we can suppose that there are no points of $G_2$ in
$B$.\\
The points of $G_1$ in $B$ can be characterized as the zero set of the 
determinant $\det dF(q,v,t)$ which is of class $C^{r-2}$.\\
We claim that
\begin{equation}\label{claim100}
\frac{\partial\det dF}{\partial t}(q,v,c(q,v))\not=0
\end{equation}
at the points of $G_1\cap B$. By the implicit function theorem, this
fact implies that $G_1\cap B$ is a regular $(n-1)$--dimensional
submanifold of class $C^{r-2}$ and we are done.

Let $\nabla$ be the covariant derivative of $(M,g)$. 
Any vector in $T_{(q,v)}UK$ can be represented by the
velocity vector $(w,u)=\nabla_s(q(s),v(s))\vert_{s=0}$ at $s=0$ of a
$C^{r-1}$ curve $(q(s),v(s))$ in $UK$, with $(q(0),v(0))=(q,v)$. Such
a curve is given by a $C^{r}$ curve $q(s)$ in $K$ with a $C^{r-1}$
unit vector field $v(s)$ defined along $q(s)$ and normal to $K$. It is
then clear that $w=q^\prime(0)$ belongs to the tangent space to $K$ at
$q$.\\
We set $u(s)=\nabla_s v(s)$ and $u(0)=u$. 
Suppose that $z(s)$ is an arbitrary vector field along $q(s)$ tangent
to $K$ with $z(0)=z$, then, by the orthogonality of $v(s)$ and $z(s)$, 
we have
\begin{equation*}
0=\frac{d}{ds}[g(z(s),\,v(s))]= g(z(s),\,u(s))+g(\nabla_sz(s)\,,v(s))
\end{equation*}
and at the point $s=0$ we get
\begin{equation*}
0=g(z,\,u)+g(\nabla_sz(s),\,v)\vert_{s=0}\,.
\end{equation*}
Introducing the shape operator $A_v:T_qK\to T_qK$ of $K$ at $q$,
relative to the unit vector $v$, we can rewrite this equation as  
\begin{equation*}
g(z,\,u)+g(A_vz,\,w)=0
\end{equation*}
and, by the symmetry property of the shape operator, 
\begin{equation*}
g(z,\,u)+g(A_vw,\,z)=0\,.
\end{equation*}
Hence, as $z$ can be chosen arbitrarily in $T_qK$, we obtain that 
$u+A_vw\in N_qK$. Notice also that, since $v$ is a unit vector,
differentiating $g(v(s),\,v(s))=1$ we obtain $g(u,\,v)=0$.\\
Resuming, the tangent space to $UK$ at the point $(q,v)$ is
represented by the pairs of vectors $(w,u)\in T_qM\times T_qM$ such
that
\begin{equation}\label{shap}
\text{$w\in T_qK$, \quad\quad $u+A_vw\in N_qK$ \quad and \quad
$g(u,\,v)=0$}\,.
\end{equation}

Consider now a vector $(w,u)\in T_{(q,v)}UK$, and the vector field
\begin{equation*}
X(t)=\partial_{UK}F(q,v,t)(w,u)
\end{equation*}
along the normal geodesic $\gamma(t)=F(q,v,t)$ with unit velocity
vector
\begin{equation*}
\gamma^\prime(t)=\partial_t F(q,v,t)\,,
\end{equation*}
where $\partial_{UK}$ denotes the partial derivative with
respect to the variable $(q,v)\in UK$.\\
The field $X$ is a Jacobi field along the geodesic $\gamma(t)$, that
is, it satisfies the following relations
\begin{equation}\label{jacobifield}
\begin{split}
\text{$X(0)=w$, \quad\quad $X^\prime(0)=u$,\quad\quad}\\
X^{\prime\prime}(t)+R(X(t),\gamma^\prime(t))\gamma^\prime(t)=0\,,
\end{split}
\end{equation}
where $R$ is the Riemann curvature operator of $M$ and we adopted the
convention of denoting with $T^\prime$ the covariant derivative along
the geodesic $\gamma(t)$ of any vector or tensor field $T$.\\
We take a basis $\left\{(w_i,u_i)\right\}$, for $i=1,\dots, n-1$,
of the tangent space $T_{(q,v)}UK$ and we construct an $n$--vector 
$\omega$ along $\gamma$ as follows,
\begin{equation*}
\omega(t)=X_1(t)\wedge\dots\wedge X_{n}(t)
\end{equation*}
where the fields $X_i(t)=\partial_{UK}F(q,v,t)(w_i,u_i)$ are Jacobi
fields and $X_n(t)=\gamma^\prime(t)$. Notice that the
relation
\begin{equation*}
X_n^{\prime\prime}(t)+R(X_n(t),\gamma^\prime(t)) \gamma^\prime(t)=0
\end{equation*}
is satisfied by the field $X_n$, since $\gamma^{\prime\prime}=0$. As
$\left\{(w_i,u_i)\right\}$ is a basis of $T_{(q,v)}UK$, proving 
equation~\eqref{claim100} is equivalent to show that
$\omega^\prime(c(q,v))\not=0$.\\
We argue by contradiction, by the Gauss Lemma (see~\cite[Chapter~3,
Section~E]{gahula}) we have that $X_n(t)$ is orthogonal to $X_i(t)$
for every $i=1, \ldots, n-1$, so we can suppose that the nonzero
vector $(w_1,u_1)$ is the generator of the kernel of
$\partial_{UK}F(q,v,c(q,v))$. Hence, $X_1(c(q,v))=0$ and $X_2(c(q,v)),
\dots, X_n(c(q,v))$ generate a subspace of dimension $(n-1)$, by the
assumption on the rank of $dF(q,v,c(q,v))$.\\
Computing $\omega^\prime(c(q,v))$ we get
\begin{equation*}
\omega^\prime(c(q,v))=X_1^\prime(c(q,v))\wedge
X_2(c(q,v))\wedge\dots\wedge X_{n}(c(q,v))
\end{equation*}
by linearity and the hypothesis that $X_1(c(q,v))=0$.\\
To conclude we need only to show that $X_1^\prime(c(q,v))$
cannot belong to the $(n-1)$--dimensional subspace generated by
$X_2(c(q,v)), \dots, X_{n}(c(q,v))$.\\
First we show that $X_1^\prime(t)$ is orthogonal to $X_n(t)$ for every 
$t$. The derivative of the function
$h(t)=g(X_1^\prime(t),\,\gamma^\prime(t))$ is zero by
$\gamma^{\prime\prime}=0$ and the last equation
in~\eqref{jacobifield}, moreover the last condition of
equation~\eqref{shap} shows that $h(0)=0$, hence the function $h$ is
identically zero.\\
Consider now the function $f(t)$ given by
\begin{equation*}
f(t)=g(X_1^\prime(t),\,X_i(t))-g(X_1(t),\,X_i^\prime(t))\,.
\end{equation*}
We have $f(0)=g(u_1,\,w_i)-g(w_1,\,u_i)$ so, 
using the second relation in~\eqref{shap} and taking into account that 
$w_i\in T_qK$ if $i\leq n-1$, we obtain 
$f(0)=g(w_1,\,A_vw_i)-g(A_vw_1,\,w_i)$ which is zero since the shape
operator is symmetric. Moreover,
\begin{align*}
f^\prime(t)=&\,g(X_1^{\prime\prime}(t),\,X_i(t))-
g(X_1(t),\,X_i^{\prime\prime}(t))\\
=&\,R(X_i(t),\gamma^\prime(t),\gamma^\prime(t),X_1(t))-
R(X_1(t),\gamma^\prime(t),\gamma^\prime(t),X_i(t))=\,0
\end{align*}
by the properties of the curvature tensor.\\
Hence, the function $f$ is identically zero and $f(c(q,v))=0$ gives
\begin{equation*}
g(X_1^\prime(c(q,v)),\,X_i(c(q,v)))=0\quad\quad\text{for $i=2,\dots, n-1$}\,,
\end{equation*}
so $X_1^\prime(c(q,v))$ is orthogonal to every vector $X_i(c(q,v))$ and
cannot belong to the subspace spanned by $\{X_i(c(q,v))\}$ for
$i=2,\dots, n$. If $X^\prime(c(q,v))$ would be zero, then by the
differential relation~\eqref{jacobifield}, we would get
$X(t)=X^\prime(t)=0$ for every $t$ and in particular for $t=0$, that
is, $(w_1,u_1)=(0,0)$ contradicting the initial hypothesis.
\end{proof}

By standard arguments of geometric measure theory
(see~\cite{fede}), the rectifiability of the cut locus has the
following immediate consequence. 

\begin{cor} The Hausdorff dimension of $\overline{Sing}$ (the cut
  locus of $K$) is at most $(n-1)$.
\end{cor}

To explain another consequence we need to introduce briefly the theory 
of functions with {\em bounded variation}, see~\cite{fede,simon} for
details. We say that a function $u:\R^n\to\R^m$ is a function with
locally bounded variation $u\in BV_{\mathrm {loc}}$, if its distributional
derivative $Du$ is a Radon measure. Such notion can be easily extended 
to maps between manifolds using smooth local charts.\\
A standard result says that the derivative of a locally
semiconcave function stays in $BV_{\mathrm {loc}}$, in view of
Proposition~\ref{carlosemicref} this implies that the vector field
$\nabla d_K$ belongs to $BV_{\mathrm {loc}}$ in the open set $M\setminus K$.

Now we define the subspace of $BV_{\mathrm {loc}}$ of functions (or vector
fields, as before) with locally {\em special bounded variation}
$SBV_{\mathrm {loc}}$ (see~\cite{amb1,amb3,amb4}).\\
The Radon measure representing the distributional  derivative $Du$ of
a function $u:\R^n\to\R^m$ with locally bounded variation can be
always uniquely separated in three mutually singular measures
\begin{equation*}
Du=D^au+Ju+Cu
\end{equation*}
where the first term is the part absolutely continuous with respect to 
the Lebesgue measure ${\mathcal L}^n$, $Ju$ is a measure concentrated
on a $(n-1)$--rectifiable set and $Cu$ (called the {\em Cantor part}) 
is a measure which does not charge the subsets of Hausdorff dimension
$(n-1)$.\\
The space $SBV_{\mathrm {loc}}$ is defined as the class of functions $u\in
BV_{\mathrm {loc}}$ such that $Cu=0$, that is, the Cantor part of the
distributional derivative of $u$ is zero.

\begin{cor} If $K$ is of class $C^r$ with $r\geq3$, the vector field
  $\nabla d_K$ belongs to the space $SBV_{\mathrm {loc}}(M\setminus K)$ of
  vector fields with locally special bounded variation.
\end{cor}

\begin{proof} Being the cut locus rectifiable, hence of Hausdorff dimension
$(n-1)$, the Cantor part of the distributional derivative of 
$\nabla d_K$ cannot be concentrated on it, so it must be concentrated
in the open set $M\setminus\left(K\cup\cutt\right)$. By point~7 of
Proposition~\ref{chiusura}, the field $\nabla d_K$ belongs to $C^{r-1}$ in
$M\setminus\left(K\cup\cutt\right)$ then, by the hypotheses, it is at
least $C^2$, hence its distributional derivative coincides with the
product of the classical derivative with the Lebesgue measure, this
shows that $Cu(M\setminus\left(K\cup\cutt\right)=0$. These two facts
together prove that $\nabla d_K$ belongs to $SBV_{\mathrm {loc}}(M\setminus K)$.
\end{proof}

A very important particular case is of our discussion is $K=\{p\}$. 
The cut locus of a point $p$ in $M$ arises naturally in various
geometric problems, its definition is due to Poincar\'e~\cite{poinc1}
and its properties were studied by many authors, see for
instance~\cite{besse1,buchner,glusin,koba1,myers1,myers2,wall,weinst1}.
A general discussion of its properties can be found in the books of
Berger~\cite{berger1} and of Gallot, Hulin,
Lafontaine~\cite{gahula}.

Because of its importance, we resume here what we got in this special
case.

\begin{prop}
Let $p$ be a point of a smooth and connected Riemannian manifold
$(M,g)$, of dimension $n$, then the squared distance function from $p$
is a locally semiconcave function on $M$ and its gradient is an $SBV$
vector field on $M$. Moreover, $d_p^2$ is $C^\infty$ in $M\setminus
Cut(p)$ which is an open neighborhood of $p$. 
The cut locus of $p$ is $C^\infty$ rectifiable with Hausdorff
dimension at most $(n-1)$ and $Cut(p)\setminus J$ is locally a finite
union of smooth hypersurfaces.
\end{prop}

\section{A Special Case}\label{analytic}

In this section we study a special example which is interesting for
the discussion of the next section. We assume that $M$ is a
two--dimensional analytic, connected and compact surface and $K$ is an 
one--dimensional embedded analytic submanifold of $M$ (a finite union
of analytic curves). As before, our analysis also applies to closed
sets $K$ with analytic boundary.\\
We look for topological results on the structure of the cut locus of
$K$ generalizing some arguments introduced to study the special case
$K=\{p\}$, see Myers~\cite{myers1,myers2}. Our goal is to show that
$\cutt$ is a finite graph and to connect its topological structure to
the differential properties of the function $d_K^2$.\\
Clearly we have that $UK$, $F$ and $F^{-1}$, when it exists, are
analytic. Notice that the fiber of $UK$ is  $U_pK \cong \{-1,1\}$.\\
The strong result given by analyticity is that the number of the
optimal focal points is finite.

\begin{lemma} The function $c:UK\to\R^+\cup\{+\infty\}$ defined in the 
  previous section is analytic in the open set where it is not $+\infty$.
\end{lemma}

\begin{proof} With the same proof of Proposition~\ref{conrect}, noticing
that when $n=2$ the set $G$ coincides with $G_1$, we can show that
the set of points $(q,v,t)\in UK\times\R^+$ where $\det dF(q,v,t)=0$
is a finite union of analytic curves. Hence, as this set is the graph
of the map $c$, we have the thesis.
\end{proof}

\begin{prop}\label{confiniti}
The set $J$ of optimal focal points of $K$ is finite.
\end{prop}

\begin{proof}
Being $M$ compact every geodesic cannot be minimal between its end
points if it is longer than the diameter of $M$, hence a minimal
geodesic joining an optimal focal point to $K$ has to be shorter than
a fixed  constant.\\ 
Consider an optimal focal point $p$ and let $F(q,v,t)$ be a 
minimal geodesic from $K$ to $p$ which has a non invertible
differential $dF(q,v,\sigma(q,v))$ (notice that in this situation we
have $c(q,v)=\sigma(q,v)$), we claim that $(q,v)$ is a critical point
of the function $c$.\\
By the Gauss Lemma (see~\cite[Chapter~3, Section~E]{gahula}), the
differential $dF(q,v,t)$ act on an element $(w,s)\in 
T_{(q,v,t)}UK\times\R^+$ as follows,
\begin{align}\label{gausslem}
dF(q,v,t)(w,s)&=\partial_{UK}F{(q,v,t)}(w)+\partial_{t}F{(q,v,t)}(s)\\
&=X+sT\nonumber
\end{align}
where the two vectors $X,T\in T_{F(q,v,t)}M$ are mutually orthogonal
and $T$ is the unit tangent vector to the geodesic $F(q,v,t)$. Taking into
account that $UK$ is locally a curve, this shows that if $dF(q,v,t)$
is singular then $\partial_{UK}F{(q,v,t)}=0$.\\
Consider now the pull--back $F^*g$ of the metric tensor
$g$ on $T_{(q,v,t)}UK\times\R^+$ via the map $F$. The set of points
$(q,v,t)$ where this form is not positive definite covers
the graph of $c$. Computing this form using 
equation~\eqref{gausslem} we have,
\begin{align}\label{arcfocal}
(F^*g)_{(q,v,t)}((w,s),(w,s))&=g_{F(q,v,t)}\left(dF(q,v,t)(w,s),\,dF(q,v,t)(w,s)\right)\nonumber\\
&=s^2+g_{F(q,v,t)}(\partial_{UK}F{(q,v,t)}(w),\,\partial_{UK}F{(q,v,t)}(w))\\
&=s^2+h(q,v,t)g_q(w,w)\nonumber
\end{align}
for a non negative function $h:UK\times\R^+\to\R$ and where $w\in
T_{(q,v)}UK$ is considered as a vector in $T_qM$. Clearly, the set of
points where the function $h:UK\times\R^+\to\R$ is equal to zero
contains the graph of $c$ by the previous discussion.\\
Suppose that $dc(q,v)\not=0$, recalling that $UK$ is a curve there
exists a small neighborhood $B$ of $(q,v)$ in $UK$ where the map $c$
is invertible and $c^{-1}$ is analytic in the open set $c(B)\subset\R^+$. 
Take a point $(r,z)\in B$ with $c(r,z)<c(q,v)$ and consider the
curve $\gamma(s):[0,c(q,v)]\to M$ defined by
\begin{equation*}
\begin{cases}
\gamma(0)=r\\
\gamma(s)=F(r,z,s)&\text{for $s\in(0,c(r,z)]$}\,,\\
\gamma(s)=F(c^{-1}(s),s)&\text{for $s\in(c(r,z),c(q,v)]$}\,,
\end{cases}
\end{equation*}
that is, after a piece of geodesic, we follow the locus of first focal
points. This is a piecewise analytic curve in $M$ starting from $K$
and ending at the point $p$. The first piece of $\gamma$ is a geodesic,
hence its length is $c(r,z)$, the second piece follows the 
locus of first focal points. Using relation~\eqref{arcfocal} we compute 
\begin{align*}
\vert\dot{\gamma}(s)\vert^2
&=g_{\gamma(s)}(dF(c^{-1}(s),s)(w(s),1),\,dF(c^{-1}(s),s)(w(s),1))\\
&=1+h(c^{-1}(s),s)g(w(s),w(s))\\
&=1
\end{align*}
for every $s\in(c(r,z),c(q,v)]$ and where
$w(s)=\frac{d(c^{-1})}{ds}(s)$. In the last equality we used the fact
that the point $(c^{-1}(s),s)$ belongs to the graph of $c$, where the
function $h$ is zero.\\
Then the length of the second piece of $\gamma$ coincides with the
variation in $s$, that is, $c(q,v)-c(r,z)$. Finally the total
length of $\gamma$ is $c(q,v)=d_K(p)$.\\
Such curve is $C^1$ since the tangent vectors of its two parts are
equal at the point $F(r,z,c(r,z))$, but it is not a geodesic for
$s\in(c(r,z),c(q,v)]$, otherwise (by uniqueness) $\gamma$ should
coincide with $F(r,z,s)$ for every $s\in(0,c(q,v)]$ and this is
impossible by construction.\\
This fact implies that there must exist a shorter
curve joining $p$ with $r$ and this is in contradiction with the
assumption that $F(q,v,t)$ is minimal, so the claim is proved.\\
Arguing by contradiction, if the set of optimal focal points 
would be infinite then in a connected component $C$  of
$\{c(q,v)<+\infty\}\subset UK$ there would be infinite points
$(q_i,v_i)\in C$ where $dc$ is zero and $F(q_i,v_i,c(q_i,v_i))=p_i$
are distinct optimal focal points. By compactness and
the initial argument on the length of a minimal geodesic from an
optimal focal point, there must exists an accumulation point of the set
$\{(q_i,v_i)\}$ in $C\subset \{c(q,v)<L\}$, for a suitable constant
$L$. Then, by the analyticity of $c$, this would imply that $dc(q,v)$
is identically zero and $c(q,v)$ is constant in the component $C$.\\
Defining a function $H:UK\to M$ by $H(q,v)=F(q,v,c(q,v))$ we have that
\begin{equation*}
{dH(q,v)}=\partial_{UK}F(q,v,c(q,v))+{\partial_t}F(q,v,c(q,v))dc(q,v)=0\,,
\end{equation*}
as $\partial_{UK}F(q,v,c(q,v))=0$ and ${dc}(q,v)=0$.\\
So the map $F(q,v,c(q,v))$ is constant in $C$. This implies
that all the points $p_i=F(q_i,v_i,c(q_i,v_i))$ coincide contradicting
the hypotheses.
\end{proof}

As the map $\sigma:UK\to M$ is continuous, the cut locus $\cutt$ is
given by a finite family of curves of kind 
$s\mapsto F(q(s),v(s),\sigma(q(s),v(s)))$ where $(q(s),v(s))$ is a
curve describing a connected component of $UK$. We say that
$p\in\cutt$ is an {\em end point} if at the point $p$ there arrives
one and only one 1--cell of points of $\cutt$. We are going to prove
that every end point is an optimal focal point.\\
First we exclude a very special case.

\begin{lemma}\label{cerchio} If at a point $p\in\cutt$ there arrive an
infinite number of minimal geodesics then all these geodesics start
from a unique connected component of $K$ which is a geodesic 
circle around $p$, that is the set of points of $M$ at a certain
distance $R$ from $p$.\\
Moreover, $p$ is an isolated point in $\cutt$, more precisely 
$\cutt \cap B_R(p)=\{p\}$. Conversely, if $p$ is an isolated point in
$\cutt$ then there is a connected component of $K$ which is a geodesic 
circle around $p$.
\end{lemma}

\begin{rem}
Notice that $p$ is an optimal focal point by  Lemma~\ref{infinite then
  focal}.\\
Since the connected components of $K$ are finite, it follows that the
isolated point of $\cutt$ are finite.
\end{rem}

\begin{proof} If $F(q_i,v_i,t)$ is the infinite family of minimal
  geodesics $F(q_i,v_i,t)$ of length $R$ ending at $p$, then all the
  distinct points $q_i$ belong to the geodesic circle of center $p$ and radius
  $R$ in $M$. The set of points $\{(q_i,v_i)\}\in UK$ clearly has an
  accumulation point, hence, by the analyticity of $UK$, the function
  $H(q,v)=F(q,v,R)$ is constantly equal to $p$ in the connected
  component of $UK$ containing such accumulation point. Again by the
  analyticity of the connected components of $UK$ and of the curves
  constituting $K$, we can conclude that the whole circle has to be a
  connected component of $K$. Hence, from every point of this circle
  there is a minimal geodesic ending at $p$ and there cannot be other
  points of $K$ inside the circle, otherwise their distance from $p$
  would be less than the radius $R$.\\
  Suppose now that $p$ is isolated in $\cutt$ and consider the open
  connected component $\Gamma$ of $M\setminus K$ which contains $p$. The 
  boundary of $\Gamma$ is a subset $K^\prime$ of $K$ and every minimal
  geodesic starting from $K^\prime$ with an initial velocity vector
  pointing toward $\Gamma$, must necessarily cease to be minimal at $p$,
  as $\cutt\cap\Gamma=\{p\}$. This last assertion follows from the fact 
  that $\Gamma$, by point~6 of Proposition~\ref{chiusura}, can be
  continuously retracted on $\cutt\cap\Gamma$, hence  $\cutt\cap\Gamma$
  is connected and then it coincides with $\{p\}$. This shows that there
  are infinite minimal geodesics from $K$ to $p$ and we can conclude as
  in the first part of the lemma.
\end{proof}

Suppose now that $p\in\cutt\setminus J$, so the number $n>1$ of
minimal geodesics $F(q_i,v_i,t)$ ending at $p$ is finite. 
Consider a small ball $B$ around $p$ in $M$, then these $n$ minimal
geodesics cut the ball $B$ in $n$ sectors that we call $S_i$. Any minimal
geodesic starting in a sufficiently small neighborhood of 
$(q_i,v_i)\in UK$ has its cut point in the ball $B$ by continuity of
the function $\sigma$, moreover this geodesic cannot cross one of the
geodesics $F(q_i,v_i,t)$ before reaching its cut point, otherwise this
latter ceases to be minimal.
Hence, considering the continuous curve of the cut points of the
geodesics starting at the points of $UK$ locally on the right side of
$(q_i,v_i)$ (remember that $UK$ is one--dimensional), we have that it 
is all contained in one of the sectors $S_i$, more precisely, by
continuity, in one of the two sectors adjacent to the geodesic $F(q_i,v_i,t)$.
This curves gives a $1$--cell of $\cutt$ approaching $p$.\\
With the same argument, considering the points locally on the left
side of $(q_i,v_i)$ we obtain another $1$--cell, in the other
sector.\\
Thus, we can conclude that the number of $1$--cells of $\cutt$ arriving 
at a point $p$ is at least the number of the sectors $S_i$, hence at
least the number of minimal geodesics from $p$ to $K$.\\
This implies that every end point of $\cutt$ where there
arrives one and only one $1$--cell, has a unique minimal geodesic to
$K$ so it has to be an optimal focal point.

Putting together these facts and Lemma~\ref{cerchio}, by
Proposition~\ref{confiniti} the end points are finite.\\
Following Myers~\cite{myers2}, this result implies that the cut locus
of $K$ is a linear graph and locally a tree, moreover the points where
the order of the graph is greater than two are finite.

Now we introduce the map $\#{\mathcal G}(p)$ from $M$ to $\NN\cup\{\infty\}$ 
{\em counting} the number of minimal geodesics from $K$ to a point $p$.

\begin{prop} An arc in $\cutt$ containing no points of $J$ and no
  interior points $p$ with $\#{\mathcal G}(p)>2$ is a regular analytic arc.
\end{prop}

\begin{proof}
Let $\gamma$ be such an arc in $\cutt$. Consider a point
$p_0\in\gamma$ with
$F(q_1,v_1,\sigma(q_1,v_1))=F(q_2,v_2,\sigma(q_2,v_2))=p_0$, by the
fact that $p_0$ is not an optimal focal point, applying the implicit
function theorem, there is an open neighborhood $B$ of $p$ in
$M\setminus K$ without optimal focal points and there exist analytic
functions $z_1,z_2:B\to UK$, $t_1,t_2:B\to\R^+$ such that,
$F(z_1(p),t_1(p))=F(z_2(p),t_2(p))=p$ for every $p\in B$ and
$z_1(B)\cap z_2(B)=\emptyset$. If $B$ is small enough, for every point 
$p$ of $B$ we have $\#{\mathcal G}(p)\leq2$, then
$t_1(p)=t_2(p)=d_K(p)$ if and only if $p\in\cutt\cap B=\gamma\cap B$.\\
The rest of the proof proceed as in Proposition~\ref{crregula}.
\end{proof}

Our last goal is to show that the {\em order} of a point $p\in\cutt$,
as a graph, is equal to $\#{\mathcal G}(p)$. The order of a point $p$
of $\cutt$ is defined as the number of distinct $1$--cells of $\cutt$
arriving at $p$. We have already seen before that that the order of
$p$ is always greater than the value $\#{\mathcal G}(p)$.\\
We now prove the opposite inequality.
 
Notice that the optimal geodesics cannot cross $\cutt$, otherwise they
cease to be minimal since they would intersect another minimal
geodesic. 
We can take a small ball $B$ around a point $p\in\cutt$ so that the
$1$--cells divide it in $n$ sectors. Let $S$ be one of these sectors
and consider a sequence of points $p_i\not\in \cutt$ all contained in $S$ and
converging to $p$ such that $F(q_i,v_i,t)$ are the minimal
geodesics relative to $p_i$. By compactness, we can suppose that the
points $(q_i,v_i)\in UK$ converge to a point $(q,v)$, hence the
minimal geodesics $F(q_i,v_i,t)$ converge to a minimal geodesic
$F(q,v,t)$ from $K$ to $p$. Being the points $p_i$ contained in $S$ 
the final part $F(q_i,v_i,t)\cap B$ of the respective minimal
geodesics have to be contained in the sector $S$ and so also the
final part $F(q,v,t)\cap B$ of the minimal geodesic for $p$.\\
Taking into account the fact that such minimal geodesic cannot
intersect the cut locus, we conclude that there is at least a minimal
geodesic for every sector $S$. Being the number of the sectors equal
to the order of $\cutt$ as a graph, we proved the opposite inequality
we claimed before.\\
Hence, there are exactly $\#{\mathcal G}(p)$ $1$--cells of the cut
locus arriving at every point $p\in\cutt$.

We summarize all the discussion of this section in the following
theorem.

\begin{teo} The set $\cutt$ is a disjoint finite union of isolated
points and linear graphs, each one locally a tree.\\
The order of every point $p\in\cutt$ equals the function $\#{\mathcal
  G}(p)$, counting the number of minimal geodesic from $p$ to $K$. In
particular, the set of points of $\cutt$ with only one or more than two
minimal geodesics is finite.\\
The set of optimal focal points in $\cutt$ is finite.\\
All the isolated points and end points of $\cutt$ are optimal focal points.\\
Considering as vertices of the graph the optimal focal points and the points
of order greater than two, the arcs connecting such vertices are
regular analytic arcs.
\end{teo}

\begin{rem}
The analysis of this section also applies with small modifications 
to the case when $K$ is a finite set of points,
considering an auxiliary set $\widetilde K$ consisting of  a family of
disjoint circles centered at the points of $K$ with a radius $R$ small
enough.
\end{rem}

\section{Singularities of Solutions of Hamilton--Jacobi Equations}

The ideas employed in the study of the distance from $K$ and
of the cut locus can be extended to analyse also the set of
singularities of viscosity solutions of general Hamilton--Jacobi
problems
\begin{equation}\label{hamjapro}
\begin{cases}
{\mathrm H}(x,du(x),u(x))=0&\text{in $\Omega\subset M$}\,,\\
u=u_0&\text{on $\partial\Omega$}\,.
\end{cases}
\end{equation}
Moreover, geometric results on the cut locus suggest conjectures 
about the viscosity solutions of these equations. We give now an
example.

Suppose that $A(x)$ is an analytic map from the closure of the unit
ball $B$ of $\R^2$ to the space of positively defined
$2\times2$--matrices.\\
We consider the following problem,
\begin{equation}\label{aequaz}
\begin{cases}
\langle A(x)\nabla u(x),\nabla u(x)\rangle=1&\text{in $B$}\,,\\
u=0&\text{on $\partial B$}\,.
\end{cases}
\end{equation}
Using arguments similar to those of Section~\ref{analytic}, it
is possible to prove that the closure of the singular set of the
viscosity solution is a finite graph.\\
If we look for the same result in the $C^\infty$ case, a
counterexample can be found as follows. It is possible to endow the
two--dimensional sphere with a $C^\infty$ metric tensor $g$ such that
the cut locus of a certain point $p$ is very wild, that is, it is not
triangulable, hence it is not a finite graph
(see~\cite{glusin}). 
Cutting away from the sphere $\SS^2$ a small geodesic disc $D$ around
$p$ whose intersection with the cut locus of $p$ is empty, and mapping
stereographically from $p$ the set $\SS\setminus D$ on $\R^2$, we have
that the closure of the singular set of the viscosity solution of 
Problem~\eqref{aequaz} in a ball of $\R^2$, where $A$ is given by the
push--forward of the metric $g$ via the stereographic projection,
coincides with a homeomorphic image of the cut locus of $p$, hence
it is not a finite graph.

However, another results on the cut locus of a point says that for a
{\em generic} $C^\infty$ metric (in the category sense, the cut locus
of every point of a surface is triangulable and has no points of order
higher than three (see~\cite{wall}). Hence, changing a little our
point of view, this discussion suggests the following conjecture.
\begin{conge} For a {\em generic} $C^\infty$ function $A(x)$ 
  from the closed unit ball $B$ of $\R^2$ to the space of positive
  definite $2\times 2$--matrices, the closure of the singular set of
  the viscosity solution of problem
\begin{equation*}
\begin{cases}
\langle A(x)\nabla u(x),\nabla u(x)\rangle=1&\text{in $B$}\,,\\
u=0&\text{on $\partial B$}
\end{cases}
\end{equation*}
is a finite graph.
\end{conge}
More in general, the same question can be asked about
Problem~\ref{hamjapro} for a {\em generic} function
${\mathrm H}$, domain $\Omega$ or boundary data $u_0$.

\bigskip

\bigskip

\noindent{\bf Acknowledgments.} The authors are grateful to Sanjoy
Mitter for having invited them to visit the Laboratory for Information
and Decision Systems at the MIT, where this work was partially done.
The authors were supported by the grant DAAL03-92-G-0115.\\
Moreover we thank Alberto Abbondandolo, Luigi Ambrosio and Michel
Delfour for many suggestions and useful discussions.

\bibliographystyle{amsplain}

\providecommand{\bysame}{\leavevmode\hbox to3em{\hrulefill}\thinspace}

\end{document}